\documentclass[12pt]{amsart}
\usepackage{amssymb}
\usepackage{amsmath}
\usepackage{graphicx}
\newtheorem{thm}{Theorem}[section]
\newtheorem{prop}{Proposition}[section]
\newtheorem{lem}{Lemma}[section]

\theoremstyle{definition}

\newtheorem{defn}{Definition}[section]
\newtheorem{rem}{Remark}[section]

\begin{document}
\title{\large Embeddability of Arrangements of Pseudocircles into the Sphere}
\author{\bf Ronald ORTNER}
\email{rortner@unileoben.ac.at}
\address{\parbox{1.4\linewidth}{Department Mathematik und Informationstechnolgie\\ Montanuniversit\"at Leoben\\
Franz-Josef-Strasse 18\\8700 Leoben, Austria}}
\begin{abstract}
An arrangement of pseudocircles is a finite set of oriented closed Jordan 
curves each two of which cross each other in exactly two points.
To describe the combinatorial structure of arrangements on closed orientable 
surfaces, in \cite{linh} so-called \textit{intersection schemes} were introduced. 
Building up on results about the latter, we first clarify the notion 
of \textit{embedding} of an arrangement. Once this is done it is shown how 
the embeddability of an arrangement depends on the embeddability of its 
subarrangements. The main result presented is that an arrangement of 
pseudocircles can be embedded into the sphere if and only if all of its 
subarrangements of four pseudocircles are embeddable into the sphere 
as well.
\end{abstract} 
\maketitle

\markboth{\sf R. Ortner}{\sf Embeddability of Arrangements of Pseudocircles into the Sphere}
\baselineskip 15pt

\section{Introduction}\label{intro}
In \cite{boko}, Bokowski used so-called hyperline sequences to introduce an 
alternative axiomatisation of oriented matroids that resulted in a new direct 
proof of the Folkman-Lawrence topological representation theorem for rank 3 
\cite{boko2} and for arbitrary rank \cite{boko3}. 
Using the hyperline sequences approach, Linhart and Ortner 
\cite{linh} introduced a generalisation of uniform oriented matroids of rank 3, 
so-called \textit{intersection schemes}, to describe the combinatorial properties 
of arrangements of pseudocircles and proved an analogue of the Folkman-Lawrence 
topological representation theorem. 
This paper continues this work by dealing with arrangements of pseudocircles
where each pair of curves intersects in exactly two crossing points. Intersection 
schemes for these particular arrangements have a special form and are introduced 
as \textit{intersection matrices}. With the aid of the latter and building up on results 
of \cite{linh} we first introduce the notion of \textit{embedding of an arrangement
in a closed orientable surface}. Then we proceed showing two minor results about
isomporphy and embeddability of arrangements in terms of their subarrangements. 
Finally, we prove that an arrangement of pseudocircles can be embedded into the sphere 
if and only if all of its subarrangements of four pseudocircles are embeddable into 
the sphere as well.

\section{Preliminaries}\label{sec:prel}
A \emph{pseudocircle} is an oriented closed Jordan curve on some closed orientable 
surface. We call a pseudocircle $\gamma$ \emph{separating}, if its complement 
consists of two connected components. With regard to one of these $\gamma$ is 
oriented counterclockwise. This component is called the \emph{interior} of 
$\gamma$, denoted by $int(\gamma)$.

\begin{defn}\label{arr}
An \emph{arrangement of pseudocircles} is a finite set of oriented closed 
Jordan curves on some closed orientable surface such that
\begin{itemize}
  \item[$(i)$] no three curves meet each other at the same point,
  \item[$(ii)$] if two pseudocircles have a point in common, they cross each other in that point,
  \item[$(iii)$] each pair of curves intersects exactly two times.
\end{itemize}
An arrangement is said to be \emph{strict} if all its pseudocircles are separating.
\end{defn}

Given an arrangement $\Gamma$ of two or more pseudocircles on an orientable 
closed surface $\mathcal S$, we may consider the intersection points of the
pseudocircles as vertices and the curves between the intersections as edges. Thus we
obtain an embedding of a graph in $\mathcal S$ which we call the 
\textit{arrangement graph}. If this induced embedding is cellular, i.e.\ all its faces
are homeomorphic to an open disc, we say that $\Gamma$ is a \textit{cellular} arrangement.

We may describe an arrangement of (labelled) pseudocircles 
$\{\gamma_1,\gamma_2,\ldots,\gamma_n\}$ as follows.
Consider a walk on each pseudocircle following its orientation beginning in an
arbitrary vertex. Whenever we meet a vertex, we note the label $i$ of the
curve $\gamma_i$ we cross provided with a sign that indicates whether $\gamma_i$
comes from the left (+) or from the right ($-$). Thus we obtain for each pseudocircle
a cyclic list of length $2(n-1)$. Ordering these lists according to the labels of the
corresponding pseudocircles one obtains an $n \times 2(n-1)$- matrix we call
the \emph{intersection matrix} of the arrangement. Generally, (abstract) intersection 
matrices may be defined as follows.

\begin{defn}\label{im}
Let ${\mathcal L} \subset \mathbb N$ be a set of $n$ labels. 
An \emph{abstract intersection matrix} is an $n \times 2(n-1)$-matrix 
whose rows are labelled with distinct elements of $\mathcal L$ in ascending order,
such that the row (with label) $i$ consists of a permutation of the elements 
$\{ +k \,|\, k\in {\mathcal L}, k\neq i \} \cup \{ -k \,|\, k\in {\mathcal L}, k\neq i \}$.

We call an intersection matrix $A$ \emph{representable} if there is an arrangement of
pseudocircles with intersection matrix $A$. If $A$ is representable with separating
pseudocircles we say that $A$ is \emph{strictly representable}. 
\end{defn}

\begin{defn}\label{cons}
An intersection matrix $A$ is \emph{consistent} if for all pairwise disjoint
$i,j,k$ the entry $\pm i$ in row $k$ is placed between $+j$ and $-j$ if and only if the
entry $\mp k$ in row $i$ is placed between $+j$ and $-j$.
\end{defn}

Linhart and Ortner \cite{linh} considered arrangements of pseudocircles of a more general 
nature: Condition (iii) of Definition \ref{arr} is relaxed such that an arbitrary finite number
of intersection points is allowed. These generalised arrangements can be described by
so-called \textit{intersection schemes}, a generalisation of our intersection matrices.
Thus, we may apply the following two main results of \cite{linh}:

\textbf{The Face Algorithm.} First, all intersection matrices are representable. For a given 
intersection matrix $A$ there is an algorithm (called the \textit{face algorithm}) that derives 
from $A$ a unique (cellular) embedding of an arrangement with intersection matrix $A$ in a 
closed orientable surface $\mathcal S_g$ of minimal genus $g$. Moreover, there can be no 
cellular embedding of an arrangement with intersection matrix $A$ in a surface $\mathcal S_{g'}$ 
of genus $g'\neq g$. Thus the combinatorial information of an arrangement is encoded in its 
intersection scheme, so that we call two arrangements of pseudocircles \textit{(combinatorially) 
isomorphic} if they can be described by the same intersection matrix. 

An arrangement $\Gamma$ is \textit{(cellularly, strictly) embeddable into} ${\mathcal S}$
if there is a (cellular, strict)  arrangement in ${\mathcal S}$ that is isomorphic to $\Gamma$.
The relation between an intersection matrix $A$ and an arrangement described by $A$ 
is analogous to that of a graph and its embedding. However, whereas a graph may be 
cellularly embeddable into several closed orientable surfaces (cf.\ \cite{gros}, p.\ 132ff), 
for an arrangement of pseudocircles there is always a unique surface it can be cellularly 
embedded into.

\begin{figure*}[h!t]
	\centering
		\scalebox{0.99}{\includegraphics{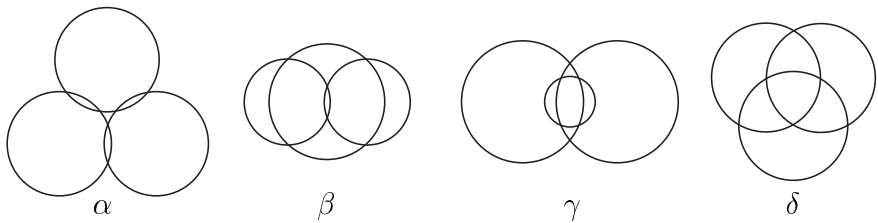}}
        \parbox{14cm}{\caption{\label{3arrs} Arrangements of three pseudocircles in the plane.}}
\end{figure*}

The face algorithm can also be used to enumerate all arrangements 
of pseudocircles. Thus Figure \ref{3arrs} shows all arrangements of three 
(counterclockwise oriented) pseudocircles in the plane (we call them 
$\alpha, \beta, \gamma, \delta$).
On closed orientable surfaces $\mathcal S_g$ there is an additional arrangement (called $\varepsilon$) 
resulting from an $\alpha$-arrangement when turning the pseudocircles inside-out, i.e.,
reversing their orientation. Actually, on the sphere and generally on $\mathcal S_g$ the arrangements 
$\alpha$, $\beta$, $\gamma$ and $\varepsilon$ can be distinguished from each other only by their
orientation. That is, one can obtain each of these arrangements from each other by reversing the
orientation of suitable pseudocircles. Unlike that, reversing any number of pseudocircles in
a $\delta$-arrangement always results in another $\delta$-arrangement. We remark that no arrangement 
of three pseudocircles can be cellularly and strictly embedded into a surface of genus $>0$.

\textbf{Characterisation of Strict Representability.} Secondly, the following theorem establishes a 
one-to-one correspondence between consistent intersection matrices and strict arrangements of 
pseudocircles.
\begin{thm}\label{char}
An arrangement is strictly embeddable into some closed orientable surface
if and only if its intersection matrix is consistent.
\end{thm}

Intersection matrices can be considered as a generalisation of oriented matroids of rank 3 
as defined via hyperline sequences by Bokowski (cf.\ \cite{boko}, p.\ 576). More precisely, 
let $A$ be an $n\times 2(n-1)$-intersection matrix. Then $A$ is a uniform oriented matroid 
of rank 3, if for all distinct $j,k$: $\pm k$ occurs in position $j$ ($1\leq j \leq n-1$) 
of a list in $A$ if and only if $\mp k$ occurs in position $j+n-1$. Given an
arrangement that can be described by a uniform oriented matroid, all its subarrangements
of three pseudocircles are of type $\delta$ (cf.\ \cite{bjoe}, p.\ 247ff).

\section{Embeddability and Isomorphy via Subarrangements}\label{sec2}
In this section we present some results concerning embeddability
and isomorphy of arrangements in terms of their subarrangements.

\begin{defn}
Let $A$ be an intersection matrix. Then \emph{submatrices of} $A$ are defined recursively
as follows. 
\begin{itemize}
	\item[(i)] The matrix arising from $A$ when deleting row $j$ and all entries 
	$\pm j$ from the other rows is a submatrix of $A$.
	\item[(ii)] If $B$ is a submatrix of $A$, then any submatrix of $B$ is also a
        submatrix of $A$.
\end{itemize}
We will shortly say \emph{$m$-submatrices} for $m\times 2(m-1)$-submatrices of $A$. 
Furthermore, $A_j^\star$ denotes the submatrix obtained from $A$ by applying rule (i).
Obviously, an $m$-submatrix of $A$ is the intersection matrix of an \emph{$m$-subarrangement} 
(subarrangement of $m$ pseudocircles) of an arrangement with intersection matrix $A$.
\end{defn}

\begin{prop}\label{isom}
Let $\Gamma,\Gamma'$ be two strict arrangements of $n\geq 4$ labelled pseudocircles
on an arbitrary closed orientable surface. Then $\Gamma,\Gamma'$ are isomorphic if and only if after a
suitable permutation of the labels they have the same set of labelled 4-subarrangements.
\end{prop}
\begin{proof}
We show that an intersection matrix with $n\geq 4$ rows is uniquely determined
by the set of its 4-submatrices. This is trivial for $n=4$. Proceeding by induction,
let $A=(a_{ik}), B=(b_{ik})$ be two intersection matrices with $n+1$ rows, such that the matrices
$A_i^\star$, $B_i^\star$ are identical for $i\in\{1,2,\ldots,n+1\}$. Now suppose that
$A$ and $B$ have different rows $j$. Since the rows are to be understood
cyclically, we may assume that $a_{j1}=b_{j1}=j+1$ (modulo $n+1$). 
Let $a_{jk}$ be the first entry in row $j$ of $A$ that is different from the 
corresponding entry $b_{jk}$ in $B$ ($k>1$). Now, because $n\geq 4$, there is an
$m\in \{1,2,\ldots,n+1\}\backslash\{j,|a_{j1}|,|a_{jk}|,|b_{jk}|\}$ (remember that $a_{j1}=b_{j1}$).
By assumption $A_m^\star=B_m^\star$ so that especially the $j$-th rows of $A_m^\star$ and
$B_m^\star$ are identical. But this leads to a contradiction:
\begin{itemize}
	\item[$\rhd$] If entries $\pm m$ neither are placed between $a_{j1}$ and $a_{jk}$ nor
        between $b_{j1}$ and $b_{jk}$, then $A_m^\star$ and $B_m^\star$ still differ in entries
        $a_{jk}^\star=a_{jk}\neq b_{jk}=b_{jk}^\star$.
	\item[$\rhd$] On the other hand, if there are entries $\pm m$ between $a_{j1}$ and
        $a_{jk}$ (or $b_{j1}$ and $b_{jk}$, respectively), then they have to appear on the
        same position in $A$ and $B$. Otherwise the $k$-th entry would not be -- as has been
        assumed -- the first one where $A$ differs from $B$. Hence, when deleting the entries
        $\pm m$ from row $j$, $A_m^\star$ and $B_m^\star$ still have different rows $j$. 
\end{itemize}
\end{proof}

\begin{rem}\label{rem:counterex2}
Proposition \ref{isom} is not true for 3-subarrangements. Figure \ref{fig:cex2}
shows two non-isomorphic arrangements (with all curves oriented counterclockwise) 
all of whose 3-subarrangements are of type $\beta$.

\begin{figure*}[h!t]
	\centering
		\scalebox{0.50}{\includegraphics{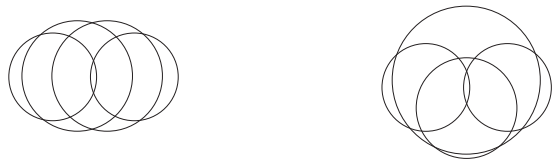}}
        \parbox{14cm}{\caption{\label{fig:cex2} An example illustrating Remark \ref{rem:counterex2}.}}
\end{figure*}

\end{rem}

As a corollary to Theorem \ref{char} one obtains the following proposition.
\begin{prop}\label{3real}
An arrangement $\Gamma$ is strictly embeddable into some closed orientable surface if and only
if all 3-subarrangements of $\Gamma$ are strictly embeddable into some closed orientable surface.
\end{prop}
\begin{proof}
By Theorem \ref{char}, it is sufficient to prove that an intersection matrix $A$ is consistent
if and only if all 3-submatrices of $A$ are consistent. Clearly, if there is an inconsistent
submatrix of $A$, $A$ itself is inconsistent. On the other hand, if $A$ is inconsistent, then
there are indices $i,j,k$, such that the entry $\pm j$ in row $k$ is placed between $+i$ and $-i$,
but the entry $\mp k$ in row $j$ is not. Hence, the 3-submatrix of $A$ consisting of the rows $i,j,k$ is
inconsistent, too.
\end{proof}

\begin{figure*}[h!t]
	\centering
		\scalebox{0.50}{\includegraphics{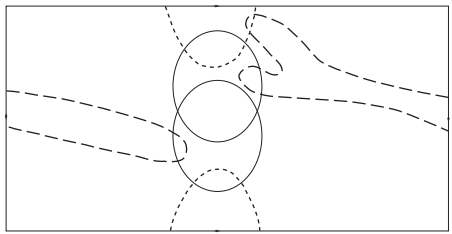}}
        \parbox{14cm}{\caption{\label{cex1} An example illustrating Remark \ref{counterex}.}}
\end{figure*}

\begin{rem}\label{counterex}
It is not sufficient that all 3-subarrangements of an arrangement $\Gamma$ are embeddable into 
some surface $\mathcal{S}$ to guarantee that $\Gamma$ is embeddable into $\mathcal{S}$ as well. 
Thus the arrangement $\Gamma$ in Figure \ref{cex1} can only be embedded into the torus (or a surface of higher 
genus), while all arrangements of three pseudocircles and hence all 3-subarrangements of $\Gamma$
are embeddable into the sphere.
\end{rem}

However, for 4-subarrangements one can prove the following theorem.
\begin{thm}\label{main}
An arrangement $\Gamma$ is embeddable into the sphere if and only if all of its 4-subarrangements
are embeddable into the sphere.
\end{thm}

\section{Proof of Theorem \ref{main}}
First, we collect some observations about embeddings of graphs and arrangements in 
surfaces of genus $>0$.

\begin{lem}\label{cycles}
Let $G=(V,E)$ be a graph that is cellularly embedded into $\mathcal{S}_g$ with genus $g>0$. 
Then the embedding contains a non-separating cycle.
\end{lem}
\begin{proof}
We apply repeatedly one of the following two operations to the embedding of $G$, until this is no longer possible:
\begin{itemize}
	\item[(i)] Remove an edge that is incident with two different faces.
	\item[(ii)] Remove a vertex of degree 1 together with the single incident edge.
\end{itemize}
Obviously, neither of the two operations has an effect on the Euler characteristic of the embedding. Thus, when it is not further possible to apply (i) or (ii), for the remaining embedded graph $G'=(V',E')$ we have $\chi(G')=\chi(G)=2-2g$. 
Since we assumed that $g>0$, $G'$ cannot be a tree. Furthermore, there are no separating cycles in $G'$, because otherwise we could apply (i). Hence, there is at least one non-separating cycle in the embedding of $G'$ and hence in
the embedding of $G$.
\end{proof}

\begin{lem}\label{sep}
Let $G=(V,E)$ be a graph (not necessarily cellularly) embedded into $\mathcal{S}_g$ ($g>0$) 
and $\mathcal{C}= (v_1,v_2\ldots,v_m)$ a non-separating cycle in $G$ $(v_i\in V,$ $v_1=v_m)$.
Furthermore, let $\mathcal{P}=(v_j,v_1',v_2',\ldots,v_\ell',v_k)$ be a path in $G$, such that
each $v_i'\in V\backslash\{v_1,\ldots,v_m\}$ and the cycle 
$\mathcal{C}_2=(v_j,v_{j+1},\ldots,v_k, v_\ell', v_{\ell-1}' \ldots,v_1',v_j)$ is separating 
(cf.\ Figure \ref{paths}).
\begin{figure*}[h]
	\begin{center}
		\scalebox{0.9}{\includegraphics{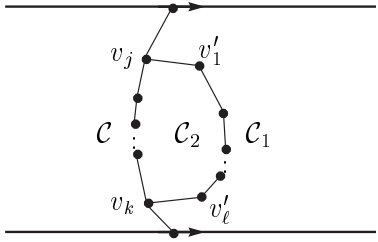}}
	\end{center}
        \caption{\label{paths}Illustration of Lemma \ref{sep}.}
\end{figure*}

Then the cycle $\mathcal{C}_1=(v_1,v_2\ldots,v_j,v_1',v_2',\ldots,v_\ell',v_k,\ldots,v_m)$
is non-separating.
\end{lem}

\begin{proof}
Assume that $\mathcal{C}_1$ is separating. Then $\mathcal{C}_1, \mathcal{C}_2$ are 
the boundary walks of two faces $F_1$ and $F_2$ whose common border corresponds 
to the path $\mathcal{P}$. It follows that $\mathcal{C}$ is the boundary walk of a 
face consisting of the union of $F_1$ and $F_2$, which contradicts our assumption 
that $\mathcal{C}$ is non-separating.
\end{proof}

\begin{lem}\label{alpha}
Let $\Gamma$ be a strict arrangement of pseudocircles cellularly embedded into 
$\mathcal{S}_g$ ($g>0$). Then (possibly after reorientation of some pseudocircles) 
$\Gamma$ contains two $\alpha$-arrangements $\Gamma_1,\Gamma_2$ such that:
\begin{itemize}
	\item[(i)] Each $\Gamma_i$ has two non-separating boundary curves, i.e.\ cycles 
		consisting of edges not contained in the interior of any pseudocircle of $\Gamma_i$.
	\item[(ii)] $\Gamma_1 \cup \Gamma_2$ can be cellularly embedded into some surface 
	of genus $>0$.
	\item[(iii)] $\Gamma_1 \cup \Gamma_2$ consists of four or five pseudocircles.
\end{itemize}
\end{lem}
\begin{proof}
According to Lemma \ref{cycles}, the arrangement graph of 
$\Gamma=\{\gamma_1,\ldots,\gamma_n\}$ contains a non-separating cycle $\mathcal C$. 
First we are going to show that we may assume that this cycle consists of edges of only 
three pseudocircles, which will give us the first $\alpha$-arrangement. Afterwards we show 
how to obtain the second one.

Given an arbitrary non-separating cycle, we may assign to it a (cyclic) sequence
of numbers $\in\{1,\ldots,n\}$, each number $i$ corresponding to a path 
consisting of consecutive edges lying on the same pseudocircle $\gamma_i$. 
We show how to obtain from an arbitrary non-separating cycle one whose 
sequence consists of only three values.  This is done in two steps. 

(1) First we delete occurring multiple entries $j$ in the sequence 
of $\mathcal{C}$. We are looking for pairs of vertices $v',v''$ of $\mathcal{C}$ on 
$\gamma_j$ such that on an oriented path $\mathcal P$ from $v'$ to $v''$ on 
$\mathcal{C}$ there are neither further intersection points with $\gamma_j$ 
nor any edges on $\gamma_j$. For each pair $v',v''$ we consider the cycle 
consisting of $\mathcal P$ together with one of the two curves connecting 
$v'$ and $v''$ on $\gamma_j$. If all these cycles were separating for all 
possible pairs $v',v''$, then some of their interiors could be composed to 
give a well-defined interior of $\mathcal C$, contradicting our assumption. 
Thus for a suitable pair $v',v''$ we find a non-separating cycle $\mathcal C'$ 
whose sequence contains only a single entry of $j$. 
Furthermore all edges of $\mathcal C'$ not on $\gamma_j$ were already 
contained in $\mathcal C$ (cf.\ Figure \ref{change1}).

\begin{figure*}[h]
	\begin{center}
		\scalebox{0.9}{\includegraphics{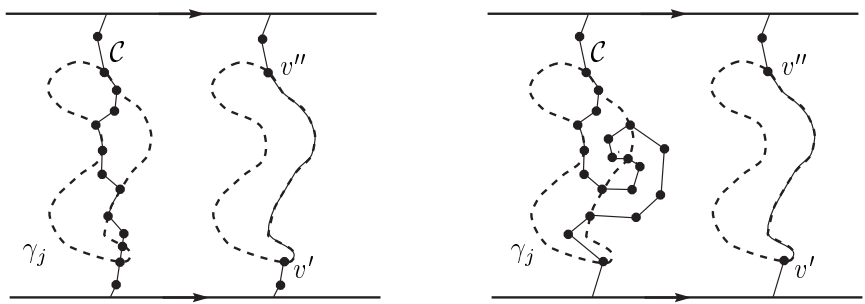}}
	\end{center}
	\caption{\label{change1}Two examples for deleting multiple entries.}
\end{figure*}

(2)  Having removed multiple occurrences in the sequence of $\mathcal{C}$, we may 
reduce it further due to the fact that each two pseudocircles intersect. Given a cycle $\mathcal{C}$ with sequence $\left\langle \ldots i, j, k,\ldots\right\rangle$ of length $\geq 4$, the basic idea is to walk on $\gamma_i$ 
ignoring the intersection with $\gamma_j$ until one arrives at an intersection with $\gamma_k$. Then continuing 
the way on $\gamma_k$ one gets back to $\mathcal{C}$ (see left picture of Figure \ref{change2}). 
There are two difficulties to consider:

First, it may happen that we cannot apply Lemma \ref{sep}, because the triangle $\Delta$ consisting 
of the detour via $\gamma_i,\gamma_k$ together with the edges on $\mathcal{C}\cap\gamma_j$ 
(i.e.\ the cycle corresponding to $\mathcal{C}_2$ in Lemma \ref{sep}) is not separating. However, in 
this case we are done, since we have found a non-separating cycle consisting only of edges on the three 
pseudocircles $\gamma_i,\gamma_j,\gamma_k$. 

\begin{figure*}[h!]
	\begin{center}
		\scalebox{0.9}{\includegraphics{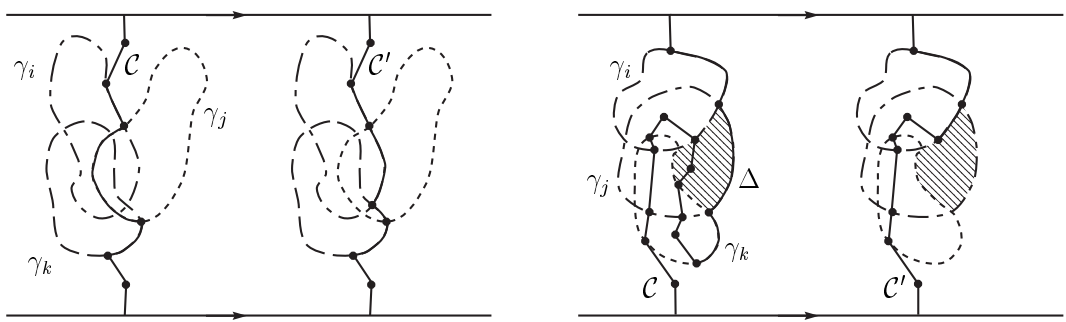}}
	\end{center}
	\caption{\label{change2}Two examples for reducing the sequence.}
\end{figure*}

Secondly, our detour may cross some part of $\mathcal{C}$ (see right picture of Figure \ref{change2}). 
Again assuming that the aformentioned triangle $\Delta$ is separating, there are at least three
cycles consisting of edges of $\Delta$ and of $\mathcal{C}$, one of which is evidently separating. 
However, it is an easy consequence of Lemma \ref{sep} that one of the other cycles has to 
be non-separating. This cycle $\mathcal{C}'$ consists only of edges of pseudocircles that also participated 
in $\mathcal{C}$. Moreover, the number of pseudocircles participating in $\mathcal{C}'$ is obviously 
smaller than in $\mathcal{C}$. Hence, repeated application of this strategy finally yields a non-separating 
cycle with sequence $\left\langle  i, j, k\right\rangle$. 

Now, all types of arrangements of three pseudocircles (cf.\ Section \ref{sec:prel}) except $\delta$ 
can be strictly embedded into surfaces of genus $>0$ so that non-separating cycles arise. 
In case $\Gamma_1=\{\gamma_i,\gamma_j,\gamma_k\}$ is an $\alpha$-arrangement it is easy
to see that this arrangement has two non-separating boundary curves $\mathcal C_1,\mathcal C_2$ 
(cf.\ Figure \ref{alpha-torus}). Otherwise, we have already mentioned in Section \ref{sec:prel} that 
such an $\alpha$-arrangement can be obtained by a suitable reorientation of some of the three 
pseudocircles $\gamma_i,\gamma_j,\gamma_k$. Thus we may assume without loss of generality
that $\Gamma_1$ is an $\alpha$-arrangement.

\begin{figure*}[h!]
	\centering
          	\scalebox{0.99}{\includegraphics{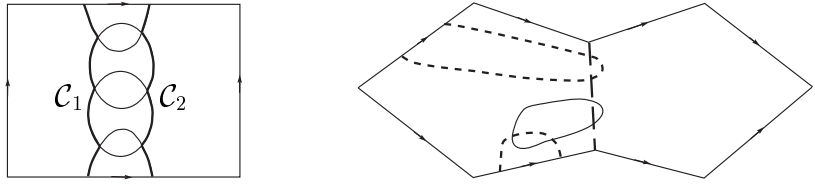}}
	\parbox{14cm}{\caption{\label{alpha-torus}Two arrangements of three pseudocircles 
embedded into some $\mathcal S_g$ ($g>0$) with non-separating cycles.}}
\end{figure*}

Obviously, in order to obtain a cellular embedding of the arrangement graph in $\mathcal{S}_g$
the cycles $\mathcal C_1$ and $\mathcal C_2$ must be connected by a simple path whose
edges are not contained in $int(\Gamma_1):=\bigcup_{\gamma\in\Gamma_1} int(\gamma)$. If 
there is a pseudocircle $\gamma_\ell\in\Gamma\setminus\Gamma_1$ that connects 
$\mathcal C_1$ and $\mathcal C_2$ with edges $\notin int(\Gamma_1)$, we have found 
(again maybe after reorientation of $\gamma_\ell$) another $\alpha$-arrangement 
$\Gamma_2$ consisting of $\gamma_\ell$ together with two pseudocircles of $\Gamma_1$. 

Otherwise, we add the pseudocircles in $\Gamma\setminus\Gamma_1$ one by one to
$\Gamma_1$ until a path as described above is established. Let $\gamma_\ell$ be the
last pseudocircle added. If $\gamma_\ell$ intersects only one of the boundary curves, 
say $\mathcal C_1$, it is clear that there must be some $\gamma_m$ that cuts $\mathcal C_2$
so that at least one intersection point of $\gamma_\ell\cap\gamma_m$ is $\notin int(\Gamma_1)$.
Then  (possibly after reorientation) $\gamma_\ell,\gamma_m$ together with an arbitrary pseudocircle 
in $\Gamma_1$ form an $\alpha$-arrangement $\Gamma_2$ embedded into $\mathcal{S}_g$ $(g>0)$ 
with two non-separating boundary curves. The case where $\gamma_\ell$ intersects both
$\mathcal C_1$ and $\mathcal C_2$ is similar. As before there is some $\gamma_m$ 
that cuts either $\mathcal C_1$ or $\mathcal C_2$ such that
$\gamma_\ell,\gamma_m$ together with an arbitrary pseudocircle $\in\Gamma_1$
forms an $\alpha$-arrangement $\Gamma_2$.

Finally, note that $\Gamma_1$ and $\Gamma_2$ were chosen such that $\Gamma_1\cup\Gamma_2$ 
can be cellularly embedded into some surface of genus $>0$.
\end{proof}

\begin{lem}\label{tech}
Let $\Gamma$ be a strict arrangement of four pseudocircles on a closed orientable surface
$\mathcal{S}_g$ of genus $g>0$ such that
\begin{itemize}
	\item[(i)] $\Gamma$ contains an $\alpha$-arrangement $\Gamma'$ with two 
			non-separating boundary curves.
	\item[(ii)] $\Gamma$ can be embedded into the sphere.
\end{itemize}
Then the single pseudocircle $\gamma\in\Gamma\backslash\Gamma'$ together with two pseudocircles
$\in\Gamma$ forms an $\alpha$-arrangement that has two non-separating boundary curves.
\end{lem}
\begin{proof}
We have to conduct an extensive case distinction. We shall see that in each single case $\gamma$ (the
dashed pseudocircle in the figures) forms together with two curves $\in\Gamma'$ (the continuous ones
in the figures) an $\alpha$-arrangement with the asked property. The figures are simplified so that
the interiors of all appearing pseudocircles are homoeomorphic to an open disc. However, the 
argumentation applies in the general case, too. Let  $\mathcal C_1, \mathcal C_2$
be the two non-separating boundary curves of $\Gamma'$.

\emph{\textbf{Case 1.}}
\emph{$\gamma$ cuts $\mathcal C_1$ and $\mathcal C_2$ such that there are edges of $\gamma$
not contained in $int(\Gamma')$ that connect a vertex $\in\mathcal C_1$ with another vertex $\in\mathcal C_2$.}\\
This can only happen if $\Gamma$ is cellularly embeddable into some surface of genus $>0$
contradicting assumption (ii).

\emph{\textbf{Case 2.}}
\emph{$\gamma$ cuts $\mathcal C_1$ and $\mathcal C_2$ such that all edges of $\gamma$
connecting a vertex $\in \mathcal C_1$ with another vertex $\in \mathcal C_2$ are contained in $int(\Gamma')$.}\\
Note that either $\mathcal C_i$ can only contain an even number of the six intersection points 
of $\gamma$ with the pseudocircles $\in\Gamma'$ so that only the following two cases may occur.

\begin{itemize}
	\item[(A)] \emph{There are two vertices of $\gamma$ on each $\mathcal C_i$:}
	First note that the two vertices on the same $\mathcal C_i$ cannot be placed
           on the same pseudocircle $\gamma_j$. Otherwise it could not be avoided that 
	$\gamma$ has more than two intersection points with $\gamma_j$. 
	\begin{figure*}[b!]
		\centering
          		\scalebox{0.75}{\includegraphics{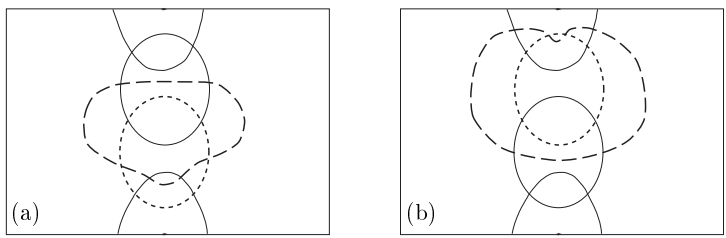}}
		\parbox{14cm}{\caption{\label{red11} Case (2.A) with all vertices of $\gamma$ on
        		$\mathcal C_1\cup\mathcal C_2$ lying on two pseudocircles.}}
	\end{figure*}
	Furthermore, since four vertices are distributed over three pseudocircles, there
	has to be a pair of vertices placed on the same pseudocircle (but on different  
	boundary curves). If all four intersection points of $\gamma$ with 
	$\mathcal C_1\cup\mathcal C_2$ lie on two pseudocircles, then the situation is as shown 
	in one of the pictures in Figure \ref{red11} (in the following we do not distinguish 
	between symmetric cases).

	Otherwise, if each pseudocircle $\in \Gamma'$ contains at least one intersection point of 
	$\gamma$ with $\mathcal C_1\cup\mathcal C_2$, then $\Gamma$ looks as in Figure \ref{red12}a.
	\begin{figure*}[b!h]
		\centering
          		\scalebox{0.75}{\includegraphics{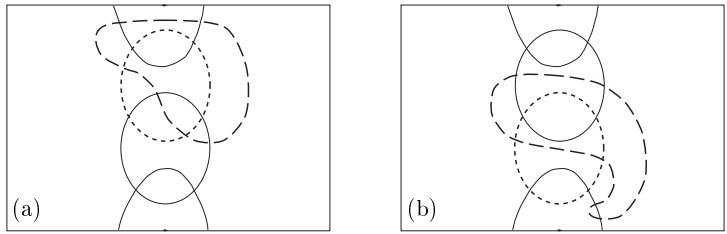}}
		\parbox{14cm}{\caption{\label{red12} Case (2.A) with vertices $\in\gamma\cap(\mathcal C_1\cup
		\mathcal C_2)$ on three pseudocircles and case (2.B).}}
	\end{figure*}
	\item[(B)] \emph{There are two vertices of $\gamma$ on $\mathcal C_1$ and 
			four on $\mathcal C_2$ (or vice versa):}
		As argued in case (A), the two vertices in $\gamma\cap\mathcal C_1$ have to be placed on different
		pseudocircles $\gamma_i,\gamma_j$. Furthermore, two of the vertices of $\gamma$ on $\mathcal 		
		C_2$ have to lie on the same pseudocircle $\gamma_k$ $(k\neq i,j)$. Hence, the other two intersection
		points of $\gamma \cap \mathcal C_2$  are placed on $\gamma_i,\gamma_j$, one on each. Thus, the
		situation is as shown in Figure \ref{red12}b.
\end{itemize}

\emph{\textbf{Case 3.}}
\emph{$\gamma$ has intersection points with $\mathcal C_1$, but not with $\mathcal C_2$ (or vice versa). }
\begin{itemize}
	\item[(A)] \emph{$\gamma$ has two intersection points with $\mathcal C_1$.}

		\begin{itemize}
		\item[(A.1)] \emph{Both of these vertices lie on the same pseudocircle $\gamma_i\in\Gamma'$:}
			The other four vertices of $\gamma$ have to be placed inside $\gamma_i$ so that
			the situation is as shown in Figure \ref{red21}a.
		\begin{figure*}[h!]
			\centering
        			\scalebox{0.75}{\includegraphics{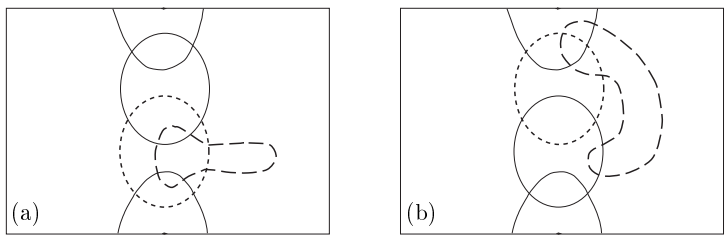}}
			\caption{\label{red21} Case (3.A.1) and case (3.B.2).}
		\end{figure*}
		\item[(A.2)] \emph{The two vertices in $\gamma\cap\mathcal C_1$ lie on different pseudocircles 									
				$\in\Gamma'$:}
			There are two possibilities dependent on how these two vertices on $\mathcal C_1$ are 		
			connected. Both are shown in Figure \ref{red22}.
		\end{itemize}
	\begin{figure*}[h!]
		\centering
		\scalebox{0.75}{\includegraphics{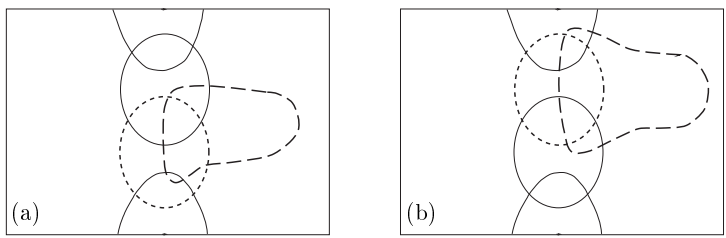}}
		\caption{\label{red22} Case (3.A.2).}
	\end{figure*}
	\item[(B)] \emph{$\gamma$ has four intersection points with $\mathcal C_1$.}
	\begin{itemize}
		\item[(B.1)] \emph{These four vertices are placed on two curves
                	$\gamma_i,\gamma_j\in \Gamma'$:} There are essentially two ways $\Gamma$
                	may look like, both shown in Figure \ref{red23}.
		\item[(B.2)] \emph{The four vertices of $\gamma$ on $\mathcal C_1$ are
                	placed on all three curves of $\Gamma'$:} This case allows only one type
                	of arrangement that can be seen in Figure \ref{red21}b.
	\end{itemize}
	\item[(C)] \emph{$\gamma$ has six intersection points with $\mathcal C_1$:} Figure 		
		\ref{red24}a shows the only possible type of arrangement satisfying this condition.
	\begin{figure*}[h!]
		\centering
  	        \scalebox{0.75}{\includegraphics{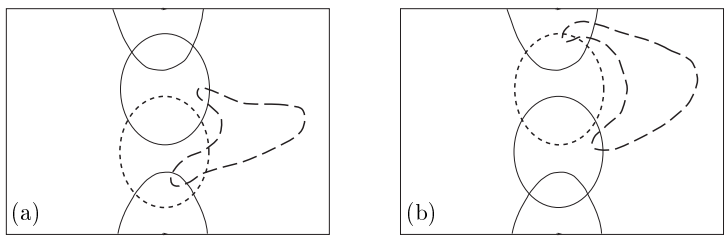}}
		\caption{\label{red23} Case (3.B.1).}
	\end{figure*}
\end{itemize}
\emph{\textbf{Case 4.}}
\emph{$\gamma$ has no intersection points with either $\mathcal C_i$.} \\
In this final case the vertices of $\gamma$ on pseudocircles $\in\Gamma'$ have to be
placed in pairs on three edges as shown in Figure \ref{red24}b.
\begin{figure*}[h!]
	\centering
          \scalebox{0.75}{\includegraphics{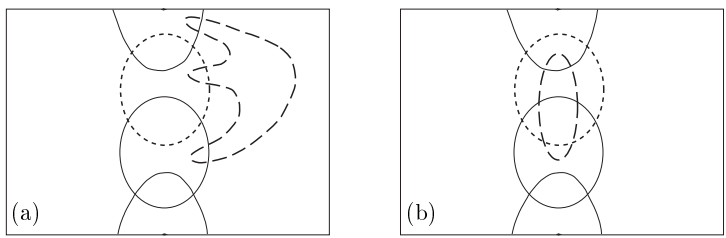}}
	\caption{\label{red24} Cases (3.C) and (4).}
\end{figure*}

In all cases we may remove the dotted pseudocircle in the figures from $\Gamma$ obtaining
an $\alpha$-arrangement with two non-separating boundary curves.
\end{proof}

\textbf{Proof of Theorem \ref{main}.}
Obviously, an arrangement can only be embeddable into the sphere if all of its 4-subarrangements
are embeddable into the sphere as well. To see that this condition is also sufficient, let $A$ be an
intersection matrix not representable on the sphere. We show that $A$ has a 4-submatrix that is 
not representable on the sphere either. If $A$ is inconsistent, we have already seen in Proposition 
\ref{3real} that there must be an inconsistent 3-submatrix and hence an inconsistent 4-submatrix 
of $A$, which is not representable in the sphere: Since pseudocircles on the sphere 
are always separating,  in this case representability and strict representability coincide. Therefore, by 
Theorem \ref{char} any inconsistent intersection matrix is not representable on the sphere.

Thus let us assume that $A$ is consistent. By Theorem \ref{char}, $A$ is strictly representable 
on some closed orientable surface $\mathcal S$ of genus $>0$. We may assume that the embedding of the 
corresponding arrangement $\Gamma$ in $\mathcal S$ is cellular. By Lemma \ref{alpha}, possibly 
after reorientation of some pseudocircles in $\Gamma$ there are two distinct $\alpha$-subarrangements
$\Gamma_1, \Gamma_2 \subseteq \Gamma$ with non-separating boundary curves such that 
$\Gamma_1\cup\Gamma_2$ consists of four or five pseudocircles and can be cellularly 
embedded into $\mathcal{S}_g$ ($g>0$). 
Now if  $|\Gamma_1 \cup \Gamma_2|=4$, we have obviously found a 4-subarrangement of $\Gamma$ 
that can be cellularly embedded into $\mathcal{S}_g$ and hence not into the sphere (cf.\ Section \ref{sec:prel}). 
Thus let us assume that $|\Gamma_1 \cup \Gamma_2| =5$, and let $\Gamma_1=\{\gamma_1, \gamma_2,\gamma_3\}$, $\Gamma_2=\{\gamma_1', \gamma_2',\gamma_3'\}$ such that 
$\gamma_1=\gamma_1'$.  Applying Lemma \ref{tech} we are going to show that one can always 
remove one of the pseudocircles from $\Gamma_1\cup\Gamma_2$ so that there still remain
two $\alpha$-arrangements that are cellularly embeddable into a surface of genus $>0$.
Thus, if $\Gamma_1 \cup \{\gamma_2'\}$ is cellularly embeddable 
into some surface of genus $>0$, we have found what we are looking for. 
Otherwise, $\Gamma_1 \cup \{\gamma_2'\}$ satisfies the conditions of Lemma 
\ref{tech} and we may infer that there are two pseudocircles 
$\gamma_{i},\gamma_{j}\in\Gamma_1$ such that 
$\Gamma'_1=\{\gamma_{i},\gamma_{j},\gamma_2'\}$ 
is an $\alpha$-arrangement with two non-separating boundary curves. For the remaining $\gamma_{k}\in\Gamma_1$ ($k\neq i,j$) there are two possibilities:
\begin{itemize}
\item[(i)] If $\gamma_{k}\neq\gamma_1$, we may remove the pseudocircle $\gamma_k$
	from $\Gamma_1\cup\Gamma_2$. This leaves the $\alpha$-arrangement $\Gamma_2$
	untouched while $\Gamma_1$ is replaced by $\Gamma_1'$, so that 
	$\Gamma_1'\cup\Gamma_2$ still is not embeddable into the sphere. 
	However, $|\Gamma_1'\cup \Gamma_2|=4$.
\item[(ii)] If $\gamma_{k}=\gamma_1$, we consider the arrangement
	$\Gamma_2\cup\{\gamma_2\}$. If it is embeddable into a surface of
	genus $>0$, we are done. Otherwise, we find analogously as described above  
	two pseudocircles $\gamma_{\ell}',\gamma_{m}'\in\Gamma_2$ such that 
        	$\Gamma_2'=\{\gamma_{\ell}',\gamma_{m}',\gamma_2\}$ is an $\alpha$-arrangement 
	whose two boundary curves are non-separating. If the remaining pseudocircle 
	$\in\Gamma_2$ is not $\gamma_1'$, the situation is as in case (i) with 
	$\Gamma_1,\Gamma_2$ interchanged.
        	Otherwise, we may remove $\gamma_1=\gamma_{1}'$ from $\Gamma_1\cup\Gamma_2$. 
	Then the $\alpha$-arrangements $\Gamma_1$, $\Gamma_2$ are replaced by
	$\Gamma_1'=\{\gamma_2,\gamma_3,\gamma_2'\}$, 
	$\Gamma_2'=\{\gamma_2',\gamma_3',\gamma_2\}$ so that $\Gamma_1'\cup\Gamma_2'$ 
	is an arrangement of four pseudocircles that still is cellularly embeddable into some surface 
	of genus $>0$.\hfill $\Box$
\end{itemize}

\section{Final Remarks}
\begin{figure*}[ht]
	\centering
		\scalebox{0.7}{\includegraphics{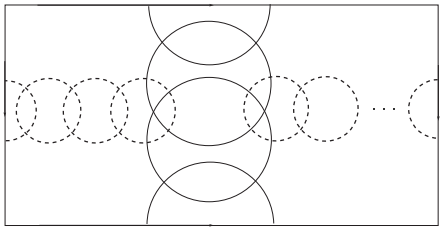}}
	{\caption{\label{cex2}A cellular generalised arrangement on the torus with all
        subarrangements embeddable into the sphere.}}
\end{figure*}

Proposition \ref{3real} can be extended to the case of connected arrangements with
relaxed condition (iii) (cf. the paragraph after Definition \ref{cons}). Unlike that, Theorem 
\ref{main} cannot be transformed into a valid version for these generalised arrangements. 
Quite to the contrary, for each $n\geq 5$ one can give a cellularly embedded generalised 
arrangement $\Gamma$ of $n$ pseudocircles on the torus, such that each ($n-1$)-subarrangement 
of $\Gamma$ is embeddable into the sphere (see Figure \ref{cex2}).

It is an interesting question whether it is possible to obtain similar results
concerning embeddability into surfaces of genus $>0$. One could e.g.\ conject that one
has embeddability into a surface $\mathcal{S}_g$ of genus $g$ if and only if all 
$(4+g)$-subarrangements are embeddable into $\mathcal{S}_g$. However, at the
moment we neither have any evidence for nor against this conjecture, and it may
be that a generalisation of Theorem \ref{main} looks totally different.

An enumeration of all 72 arrangements of four pseudocircles that can be embedded into the 
sphere can be found in \cite{ortn}. By the way, each of these arrangements can be realised 
with proper circles (cf.\ Appendix B of \cite{ortn}).

\end{document}